\documentclass[a4paper,10pt]{article}
\usepackage{amsmath,amsthm,amssymb,latexsym}
\usepackage[latin1]{inputenc}   

\usepackage{graphicx}

\newcommand{\smallbox}{\rule{5pt}{5pt}}

\newtheorem{corollary}{Corollary}[section]
\newtheorem{conjecture}{Conjecture}[section]

\newtheorem{lemma}{Lemma}[section]

\newtheorem{thm}{Theorem}[section]
\newtheorem{theorem}[thm]{Theorem}

\theoremstyle{remark}
\newtheorem*{remark}{Remark}

\def\ignore#1{}
\def\feliu#1{}

\newcommand{\N}{\mathbb{N}}

\newcommand{\hyp}[1]{Q_{#1}}
\newcommand{\gr}[1]{\mbox{\rm Gr}_{#1}}
\newcommand{\edg}[2]{#1#2}
\newcommand{\idim}{$i$-th dimension}

\begin{document}
\title{Parity balance of the \idim\  edges in Hamiltonian cycles of the hypercube}

\author{Feli\'u Sagols \footnote{ Mathematics Department, CINVESTAV-IPN
Mexico City. fsagols@math.cinvestav.mx}, 
Guillermo Morales-Luna \footnote{Computer Science Department, CINVESTAV-IPN 
Mexico City. gmorales@cs.cinvestav.mx}}

\date{\today}

\maketitle

\thispagestyle{empty}

\begin{abstract}
Let $n\geq 2$ be an integer, and let $i\in\{0,\ldots,n-1\}$. An \idim\  edge in the $n$-dimensional hypercube $\hyp{n}$ is an edge $\edg{v_1}{v_2}$ such that $v_1,v_2$ differ just at their $i$-th entries. The parity of an \idim\  edge $\edg{v_1}{v_2}$ is the number of 1's modulus 2 of any of its vertex ignoring the $i$-th entry. We prove that the number of \idim\  edges appearing in a given Hamiltonian cycle of $\hyp{n}$ with parity zero coincides with the number of edges with parity one. As an application of this result it is introduced and explored the conjecture of the inscribed squares in Hamiltonian cycles of the hypercube: Any Hamiltonian cycle in $\hyp{n}$ contains two opposite edges in a $4$-cycle. We prove this conjecture for $n \le 7$, and for any Hamiltonian cycle containing more than $2^{n-2}$ edges in the same dimension. This bound is finally improved considering the equi-independence number of $\hyp{n-1}$, which is a concept introduced in this paper for bipartite graphs.

\end{abstract}

\section{Introduction} \label{s:introduction}

Let $n \ge 1$ be an integer number. The {\em $n$-dimensional hypercube}, denoted $\hyp{n}$, is the graph having $\{0,1\}^n$ as set of vertexes, two of them being joined by an edge if they differ just in one of their entries. It is well known that hypercubes are all Hamiltonian graphs. The canonical examples of Hamiltonian cycles in the hypercube are the so called {\em binary Gray codes}~\cite{Sa97}.

The study of structural properties of the Hamiltonian cycles of the hypercube have allowed the solution of relevant problems. For instance, the proof of Kreweras's conjecture by Fink in 2007 (see~\cite{Fink07,Fink09}) about the extensibility of any perfect matching on the hypercube allowed Feder and Subi in 2009 to find tight bounds on the number of different Hamiltonian cycles in the hypercube~\cite{Feder09}; the same structural property allowed Gregor in 2009 to prove that perfect matchings on subcubes can be extended to Hamiltonian cycles of hypercubes~\cite{Gr09}. Something similar has been done with Hamiltonian paths, for example Chen proved in 2006 (see~\cite{chen06}) that any path of length at most $2n-1$ is a sub-path of a Hamiltonian path in $\hyp{n}$ and he used this basic structural property to prove that, for $n \ge 3$, $\hyp{n}$ is $(2n-3)$-path bipancyclic but it is not $(2n-2)$-path bipancyclic. Moreover, he proved that a path $P$ of length $k$, with $2 \le k \le 2n-3$, lies in a cycle of length $(2k-2)$ if and only if $P$ contains two edges in the same dimension~\cite{chen06}. 

Our objective here is to extend the structural knowledge about Hamiltonian cycles of the hypercube by classifying permissible sets of edges in Hamiltonian cycles, in terms of the dimension to which these edges belong and their parities (see Theorem~\ref{p:bipartition}).

As an application we conjecture (see Conjecture~\ref{c:squaresConjecture}) that any Hamiltonian cycle in the hypercube contains the opposite edges of a $4$-cycle, and we prove some particular instances of this conjecture (see Theorem~\ref{p:bipartition}). In the 80's Erd\"os~\cite{Er84} conjectured that any subgraph of $\hyp{n}$ asymptotically containing at most half edges of the whole $\hyp{n}$ is $4$-cycle free, recent advances on this conjecture appeared in~\cite{AxM06} and~\cite{Th09}. The substantial difference of Erd\"os conjecture and ours is that Erd\"os conjecture poses conditions on the maximum number of edges in a $4$-cycle-free graph and our conjecture is about the existence of $4$-cycles with opposite edges in any Hamiltonian cycle of the hypercube.

As a first approach to explore Conjecture~\ref{c:squaresConjecture} it is introduced the notion of \idim\  graph which in turn is isomorphic to $\hyp{n-1}$ and allows us to translate our  decision problem (whether there is a $4$-cycle within any Hamiltonian cycle $h$ in the hypercube) into verifying that the maximum number of edges in a dimension in $E(h)$ is greater than the independence number of $\hyp{n-1}$ (see Theorem~\ref{p:independence}). Then the last bound is improved as the equi-independence number of $\hyp{n-1}$ (see Corollary~\ref{c:equiIndependence}). The equi-independence number of a bipartite graph $G$ is the cardinality of the maximum independent set in $G$ containing the same number of elements in each bipartition class; we prove in Theorem~\ref{p:equiIndBipartite} that finding the equi-independence number of a graph is polynomial time reducible to the independence number computation problem.

The outline of the paper is the following: In Section~\ref{s:preliminaries} we recall very basic notions of Graph Theory and we introduce some notation, in Section~\ref{s:coloring} we introduce and study \idim\  graphs. In Section~\ref{s:balance} we introduce the notion of {\em chromatic vector} and some other operators defined over Hamiltonian cycles in the hypercube, and we prove Theorem~\ref{p:bipartition} which is the main result in this paper. In Section~\ref{s:squares} we state and discuss the inscribed square conjecture, we introduce the equi-independence concept for bipartite graphs and we report our final progress on this conjecture. In the conclusions, we suggest additional applications to Theorem~\ref{p:bipartition} and we pose a list of open conjectures and problems.

\section{Preliminaries} \label{s:preliminaries}

Let $G=(V,E)$ be graph with set of nodes $V$ and set of edges $E$. Let us recall the following elementary notions. A {\em $k$-cycle} in $G$ is a sequence of pairwise different nodes $v_0v_1\ldots v_{k-1}$ such that for each index $i$, $\edg{v_i}{v_{i+1}}$ is an edge in $E$ (index addition is taken modulus $k$). A {\em Hamiltonian cycle} is a $k$-cycle, where $k=\mbox{card}(V)$ is the order of the graph. The graph $G$ is {\em Hamiltonian} if it possesses a Hamiltonian cycle. A non-empty set $I\subset V$ is {\em independent} if no pair of different elements in $I$ is an edge: $\forall u,v\in I$ $\left[u\not=v\ \Rightarrow\ \edg{u}{v}\not\in E\right]$. A {\em maximal independent set} is an independent set which is maximal with respect to set-inclusion. The {\em independence number} $\alpha(G)$ of $G$ is the number of vertexes in a largest independent set:  $\alpha(G)=\max\{\nu|\ \exists I\subset V:\ I\mbox{ independent }\ \&\ \mbox{card}(I)=\nu\}$. Any independent set consisting of  $\alpha(G)$ vertexes is called a {\em maximum independent set}.

Let $n \ge 2$ be an integer. It may be assumed that $E(\hyp{n})$ consists of pairs $\edg{v}{(v+e_i)}$, where $e_i$ is the $i$-th canonical basic vector. Clearly, $\edg{v}{(v+e_i)}=\edg{u}{(u+e_j)}$ if and only if $i=j$ and either $u=v$ or $u=v+e_i$ (addition is integer addition modulus 2).

Any 4-cycle in $\hyp{n}$ has thus the form $v(v+e_i)(v+e_i+e_j)(v+e_j)$, with $v\in V(\hyp{n})$ and $0\leq i<j\leq n-1$.

The order of $\hyp{n}$ is thus $2^n$ and its independence number is $\alpha(\hyp{n})=2^{n-1}$.

The hypercube $\hyp{n}$ is Hamiltonian and the {\em binary Gray code}~\cite{Sa97}, $\gr{n}$, is a Hamiltonian cycle. As a sequence, this code is determined recursively by the following recurrence:
$$\gr{1} = [0,1]\ \ ,\ \ \gr{n} = \mbox{join}(0*\gr{n-1},1*\mbox{rev}(\gr{n-1}))$$
(join and rev are respectively list concatenation and list reversing, $*$ is a {\em prepend} map: $b*\mbox{list}$ prepends the bit $b$ to each entry at the list). In general, we will follow the notions and notations in Diestel's textbook~\cite{Die05}.

\section{Graphs along a dimension}\label{s:coloring}

Let $n\geq 2$ be an integer. Let $i$ be an integer in $\{0, \ldots, n-1\}$. An {\em \idim\ } edge in $\hyp{n}$ is and edge of the form $u(u+e_i)$ where $u \in V(\hyp{n})$. 

For each $i\in\{0,\ldots,n-1\}$ let ${\cal D}_{in}$ be the $i$-{\em dimension graph} whose vertexes are the \idim\  edges of $\hyp{n}$, and whose edges are the pairs of \idim\  edges forming a 4-cycle within $\hyp{n}$:
\begin{equation}
\edg{\edg{v}{(v+e_i)}}{\edg{u}{(u+e_i)}}\mbox{ is an edge in }{\cal D}_{in} \  \Longleftrightarrow \ 
\exists j\not=i:\ \edg{u}{(u+e_i)} = \edg{(v+e_j)}{(v+e_j+e_i)}. \label{eq.bied}
\end{equation}

\begin{theorem} \label{p:isomorphism}
Let $n \ge 2$ be an integer. For each $i\in\{0,\ldots,n-1\}$ the graph ${\cal D}_{in}$ is isomorphic to $\hyp{n-1}$.
\end{theorem}

\noindent {\em Proof:} The map $f_{in} : V({\cal D}_{in}) \rightarrow V(\hyp{n-1})$, $v(v+e_i)\mapsto f_{in}(v(v+e_i)) = \rho_{in}(v)$, where $\rho_{in}$ is deletion of the $i$-th entry, is a graph isomorphism. \hfill
$\smallbox$
\medskip

Hence, being all isomorphic to $\hyp{n-1}$, the \idim\  graphs ${\cal D}_{in}$ are pairwise isomorphic and each ${\cal D}_{in}$ is bipartite. 

Let $\mbox{\it par}_n:\hyp{n}\to\{0,1\}$ be the map such that $\mbox{\it par}_n(v)$ is the parity of the Hamming weight (number of entries equal to 1) of the vector $v\in\hyp{n}$, and for any index $i\in\{0,\ldots,n-1\}$ let $\mbox{\it par}_{in} = \mbox{\it par}_{n-1}\circ\rho_{in}$ where $\rho_{in}$ is the map that suppresses the $i$-th entry at any $n$-dimensional vector. In other words, $\mbox{\it par}_{in}(v)$ is the parity of the vector resulting by suppressing the $i$-th entry in $v$.
We have that the bipartition classes of ${\cal D}_{in}$ are realized as:
$$\edg{v_0}{v_1}\,,\,\edg{v_2}{v_3}\mbox{ in the same bipartition class }\Longleftrightarrow\ \mbox{\it par}_{n-1}\circ f_{in}(\edg{v_0}{v_1}) = \mbox{\it par}_{n-1}\circ f_{in}(\edg{v_2}{v_3}). $$
Or equivalently,
\begin{equation}
\edg{v}{(v+e_i)}\,,\,\edg{u}{(u+e_i)}\mbox{ in the same bipartition class }\Longleftrightarrow\ \mbox{\it par}_{in}(v) = \mbox{\it par}_{in}(u). \label{eq.bipa}
\end{equation}
According to the common value at the right side of relation~(\ref{eq.bipa}) we will refer to the partition classes as $0$- and $1$-{\em bipartition classes} of ${\cal D}_{in}$. 

Let $h\,=\,h_0\cdots h_{2^n-1}$ be a Hamiltonian cycle of $\hyp{n}$. For each $\iota\in\{0,\ldots,2^n-1\}$ let us associate the dimension $i$ of the edge $\edg{h_\iota}{h_{\iota+1}}$ as a color of the starting vertex (in the sense of the cycle) ${h_\iota}$. Let us denote $\chi_h:V(\hyp{n})\to\{0,\ldots,n-1\}$ this vertex coloring induced by the Hamiltonian cycle $h$.

Let $c(h)=\left(c_i\right)_{i=0}^{n-1}\in\N^n$ be the vector such that for each $i$, $c_i$ is the number of vertexes colored $i$ by the Hamiltonian cycle $h$. The vector $c(h)$ is called the {\em chromatic vector} of $h$.

$c:\{\mbox{Hamiltonian cycles}\}\to\N^n$ is a well defined map. A full characterization of the image $C\subset\N^n$ of this map is out of the scope of this paper. However, some necessary conditions for chromatic vectors are asserted at the following

\begin{lemma} \label{l:basicChromatic}
Let $n\geq 2$ be an integer and let $c\in C$ be the chromatic vector of some Hamiltonian cycle $h$, $c=c(h)$. Then:
\begin{enumerate}
 \item All entries of $c$ are even integers.
 \item No entry at $c$ is zero.
 \item The sum of the entries of $c$ is equal to $2^n$.
 \item The greatest entry in $c$ is lower or equal than $2^{n-1}$.
 \item The lowest entry in $c$ is greater or equal than $2$.
 \item Any permutation of $c$ is the chromatic vector of a Hamiltonian cycle of $\hyp{n}$.
\end{enumerate}
\end{lemma}

\noindent {\it Proof:}  The first assertion follows from the fact that on the hypercube $\hyp{n}$, $\sum_{v\in\hyp{n}} v = 0\in V(\hyp{n})$ (sum is addition modulus 2). The second because if $c_i=0$, for some $i$, then the nodes at the Hamiltonian $h$ either all lie at the semi-space $v_i=0$ or at the opposite semi-space $v_i=1$, which is not possible. Third assertion holds because each vertex has assigned a color. Fourth assertion follows from the fact that no consecutive pair of edges in $h$ can be parallel, or, in other words, if $h_jh_{j+1}$ has the same direction as $h_{j+1}h_{j+2}$, for some $j \in \{0, \ldots, 2^{n}-1\}$, then necessarily $h_j = h_{j+2}$ but that is impossible. Fifth assertion is a consequence of the first and the second. For the last assertion, let $\pi:\{0,\ldots,n-1\}\to\{0,\ldots,n-1\}$ be a permutation, then $c(\pi\circ h)=\pi(c(h))$ and $\pi\circ h$ is a Hamiltonian cycle. \hfill$\smallbox$

Let $h\,=\,h_0\cdots h_{2^n-1}$ be a Hamiltonian cycle in the hypercube and let $i\in\{0,\ldots,n-1\}$. By rotating and reverting, if necessary, the indexes in $\{0,\ldots,2^n-1\}$, we may assume that $i$ is the color assigned to $h_0$ and that ${h_0}_i = 0$ (indeed, $\edg{h_0}{h_1}$ is an \idim\  edge). Let us define the following operators:
\begin{description}
\item[inverse image index list. ]  $\lambda_i(h)=[j_0,j_1,\ldots,j_{c_i-1}]$ ($j_0=0$):  the list of indexes colored $i$ at the Hamiltonian $h$.  This is the list of indexes $\ell\in\{0,\ldots,2^n-1\}$ such that $\chi_h(h_\ell)=i$. In other words, it is the list of indexes in $h$ where an \idim\  edge starts,
\item[inverse image lists. ]   $\eta_i(h)=[h_{j_0},h_{j_1},\ldots,h_{j_{c_i-1}}]$, and 
$\xi_i(h)=[\edg{h_{j_0}}{h_{j_0+1}},\edg{h_{j_1}}{h_{j_1+1}},\ldots,\edg{h_{j_{c_i-1}}}{h_{j_{c_i-1}+1}}]$. The list $\eta_i(h)$ consists of the $i$-th colored vertexes at $\hyp{n}$, and the list $\xi_i(h)$ consists of the \idim\  edges at $\hyp{n}$.\end{description}
Let us study the maps $\lambda_i$, $\eta_i$, $\xi_i$.

\section{Parity balance}\label{s:balance}

Let $n\geq 2$ be an integer.  The $n$-dimensional hypercube $\hyp{n}$ is bipartite and its bipartite classes are the collections of vertexes of odd and even parity: $P_{0n} = \mbox{\it par}_n^{-1}(0)$, $P_{1n} = \mbox{\it par}_n^{-1}(1)$. Let $h\,=\,h_0\cdots h_{2^n-1}$ be a Hamiltonian cycle in the $n$-dimensional hypercube $\hyp{n}$.

\begin{remark}
The sequence $\mbox{\it par}_n\circ h$ consists of alternating values 0 and 1.
\end{remark}

\begin{theorem} \label{p:bipartition}
 Let $n \ge 2$ be an integer, let $i\in\{0, \ldots, n-1\}$ and let $h$ be a Hamiltonian cycle in $\hyp{n}$. Then half of the \idim\  edges appearing in $h$, namely at the list $\xi_i(h)$, lie in the $0$-bipartition class of ${\cal D}_{in}$ and the other half in the $1$-bipartition class. 
\end{theorem}

\noindent {\em Proof:} The first statement at Lemma~\ref{l:basicChromatic} asserts that $c_i = |\xi_i(h)|$ is indeed an even integer. Without loss of generality we may assume that $\edg{h_0}{h_1}$ is an \idim\  edge and the $i$-th entry of the starting vertex $h_0$ has value $0$, i.e. ${h_0}_i = 0$.

For a bit $b \in \{0,1\}$, let $\hyp{bin}$ be the subgraph of $\hyp{n}$ induced by those vertexes with value $b$ at their $i$-th entry. Both $\hyp{0in}$ and $\hyp{1in}$ are isomorphic to $\hyp{n-1}$ and they are bipartite as well.

Let $\lambda_i(h)=[j_0,j_1,\ldots,j_{c_i-1}]$ ($j_0=0$) be the index list corresponding to $i$-colored nodes at $h$. From relation~(\ref{eq.bipa}) we have
\begin{equation}
\forall k<c_i,b \in \{0,1\}\ \left[\mbox{\it par}_{in}(h_{j_k}) = \mbox{\it par}_{in}(h_{j_k+1}) = b\ \Longleftrightarrow\ \edg{h_{j_k}}{h_{j_k+1}}\in b\mbox{-bipartition class}\right]. \label{eq.bipao}
\end{equation}
Let us consider the parity list $\left[\mbox{\it par}_{in}(h_{j_k})\right]_{k=0}^{c_i-1}$. By relation~(\ref{eq.bipao}), the theorem will be proved by showing that half values in this list are 0 and half are 1.

\begin{figure}
     \centerline{\includegraphics[height=8cm]{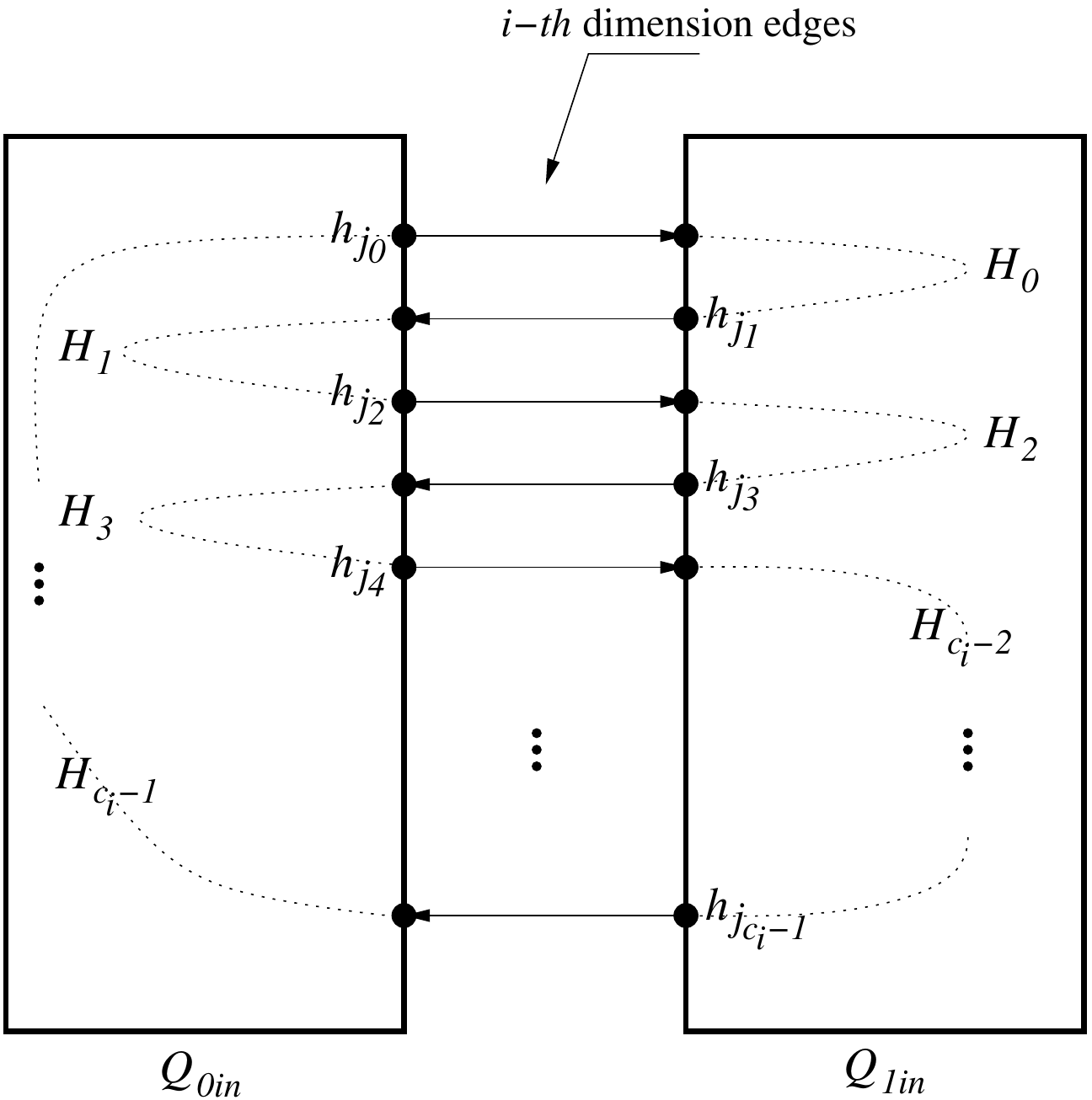}}
\caption{The traces of a Hamiltonian cycle on the parallel semi-hypercubes with constant $i$-th entry. \label{fg.01}}
\end{figure}

Since ${h_0}_i = 0$, for any even $k\in\{0, 1, \ldots, c_i-1\}$ the path $H_k\,=\,h_{j_k+1} \cdots h_{j_{k+1}}$ lies entirely in $\hyp{1in}$. Let $V(H_k)$ be the collection of vertexes appearing in $H_k$. Indeed, 
\begin{equation}
\left\{V(H_{2k_1})\right\}_{k_1=0}^{\frac{c_i}{2}-1} \mbox{ is a partition of } V(\hyp{1in}). \label{eq.partitionEven}
\end{equation}
Similarly, 
\begin{equation}
\left\{V(H_{2k_1+1})\right\}_{k_1=0}^{\frac{c_i}{2}-1} \mbox{ is a partition of } V(\hyp{0in}). \label{eq.partitionOdd}
\end{equation}
Figure~\ref{fg.01} displays a diagram of this situation.
According to relation~(\ref{eq.bipao}), the following equivalences hold:
\begin{eqnarray}
\mbox{\it par}_{in}(h_{j_k+1}) \not = \mbox{\it par}_{in}(h_{j_{k+1}}) &\Longleftrightarrow& \mbox{the length of $H_k$ is odd,} \label{eq.051}\\
\mbox{\it par}_{in}(h_{j_k+1}) = \mbox{\it par}_{in}(h_{j_{k+1}}) &\Longleftrightarrow& \mbox{the length of $H_k$ is even.}\label{eq.052}
\end{eqnarray}

If the inequation at the left side of~(\ref{eq.051}) holds then the endpoints of $H_k$ have different parities and they does not introduce any imbalance on the number of zeros and ones in 
$\left[\mbox{\it par}_{in}(h_{j_k})\right]_{k=0}^{c_i-1}$ (recall that by definition $\mbox{\it par}_{in}(h_{j_k+1}) = \mbox{\it par}_{in}(h_{j_k})$). 
If the equation at the left side of~(\ref{eq.052}) holds with value $b\in\{0,1\}$ for some even number $k=k_1$, then $h_{k_1}$ has an odd number of vertexes and in consequence one vertex less in the $\overline{b}$-bipartition class of $\hyp{1in}$ than in the $b$-bipartition class. Since relation~(\ref{eq.partitionEven}) is true this deficit must be compensated with the existence of some even value $k=k_2$ such that $\mbox{\it par}_{in}(h_{j_{k_2}+1}) = \mbox{\it par}_{in}(h_{j_{k_2+1}}) = \overline{b}$. In other words
\begin{eqnarray*}
\mbox{card}\left(\{ k \in \{0, 2, \ldots, c_i-2\} | \mbox{\it par}_{in}(h_{j_k+1}) = 0 = \mbox{\it par}_{in}(h_{j_{k+1}}\}\right) &=& \\ 
\mbox{card}\left(\{ k \in \{0, 2, \ldots, c_i-2\} | \mbox{\it par}_{in}(h_{j_k+1}) = 1 = \mbox{\it par}_{in}(h_{j_{k+1}}\}\right) &&
\end{eqnarray*}
Consequently $\left[\mbox{\it par}_{in}(h{j_k})\right]_{k=0}^{c_i-1}$ shall have the same number of 0's and 1's. \hfill $\smallbox$

The above result entails some necessary conditions for a vector to be of the form $c(h)$ for some Hamiltonian cycle $h$.

For a Hamiltonian cycle $h$ on $\hyp{n}$ and an index $i\in\{0, \ldots, n-1\}$ let us define:
\begin{description}
\item[chromatic-segments list. ]  $\delta_i(h) = [j_1-j_0,j_2-j_1,\ldots,2^n-j_{c_i-1}]$: the difference among the right shift of $\lambda(h)$ and $\lambda(h)$ itself,
\item[parity list. ]  $\beta_i(h)=[b_0,b_1,\ldots,b_{c_i-1}]$: bit list defined recursively as follows:
$$b_0 = \mbox{\it par}_{in}(h_0)\ \ \&\ \ \forall k\in\{1,\ldots,c_i-1\}:\ b_k = (b_{k-1} + j_k - j_{k-1}+1) \bmod 2.$$
\end{description}

\begin{corollary} \label{c:parityVectors}
Let $n \ge 2$ be an integer and $i\in\{0, \ldots, n-1\}$, and let $h$ be a Hamiltonian cycle on $\hyp{n}$. Then the $i$-th parity list of $h$, $\beta_i(h)$, is a balanced list, i.e. the number of $0$ and $1$ values coincide.
\end{corollary}

\noindent {\em Proof:}
A simple induction proves that $\beta_i(h) = \left[\mbox{\it par}_{in}(h_{j_k})\right]_{k=0}^{c_i-1}$. The conclusion follows from the proof of Theorem~\ref{p:bipartition}. \hfill $\smallbox$

\begin{corollary}\label{l:chromaticsegmets}
Let $n \ge 2$ be an integer and $i\in\{0, \ldots, n-1\}$. The addition of the even numbered entries at the chromatic-segments list $\delta_i(h)$ coincides with the addition of its odd numbered entries and both are equal to $2^{n-1}$:
\begin{equation}
\sum_{k_1=0}^{\frac{c_i}{2}-1}(\delta_i(h))_{2k_1} = 2^{n-1} = \sum_{k_1=0}^{\frac{c_i}{2}-1}(\delta_i(h))_{2k_1+1} \label{eq.igual}
\end{equation}
\end{corollary}

\noindent {\em Proof:} For any $k \in \{0, \ldots, c_i-1\}$ the path $H_k$ in the proof of Theorem~\ref{p:bipartition} has $\delta_k(h)$ vertexes. The first (resp. last) expression at relation~(\ref{eq.igual}) corresponds to the sum of the lengths of paths $h_k$ with even (resp. odd) index and the result follows from relations~(\ref{eq.partitionEven}) and~(\ref{eq.partitionOdd}). \hfill $\smallbox$

\section{Inscribed squares in Hamiltonian cycles} \label{s:squares}

We introduce here an application of Theorem~\ref{p:bipartition}. We will establish sufficient conditions over Hamiltonian cycles, in terms of their chromatic vectors, to guarantee  that there are 4-cycles within the Hamiltonian cycles.

\begin{figure}
\begin{center}
\includegraphics[height=5cm]{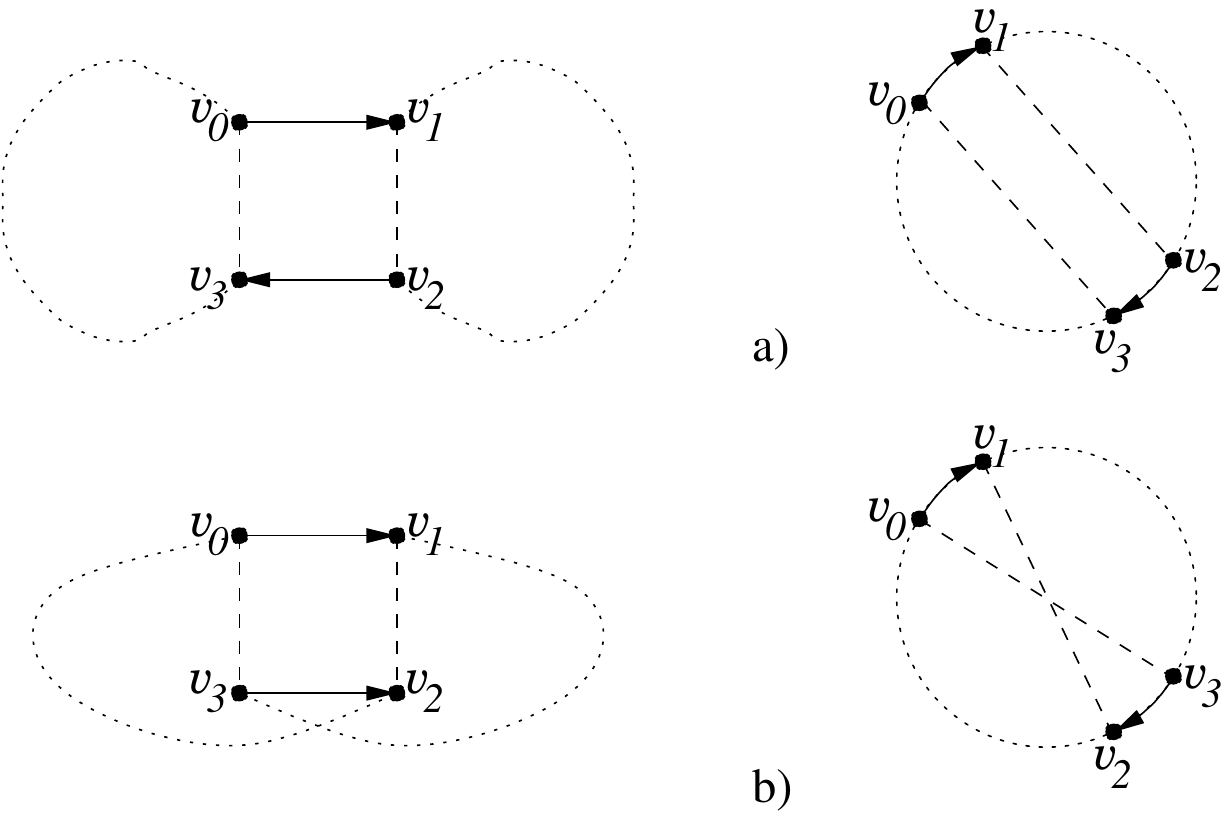}
\end{center}
\caption{Hamiltonian cycles possessing squares. (a) Straight squares. (b) Twisted squares. \label{fg.02}}
\end{figure}

Let $n\geq 2$ be an integer and let $h\,=\,h_0\cdots h_{2^n-1}$ be a Hamiltonian cycle of $\hyp{n}$. Let $\edg{v_0 v_1}{v_2 v_3}$ be a 4-cycle in the hypercube $\hyp{n}$. Let us introduce the following definitions:
\begin{enumerate}
\item $\edg{v_0 v_1}{v_2 v_3}$ is a {\em straight square} within $h$ if there exist indexes $i,j \in \{0, \ldots, 2^n-1\}$ such that 
$$h_i h_{i+1} h_j h_{j+1} = v_0v_1v_2v_3.$$
\item $\edg{v_0 v_1}{v_2 v_3}$ is a {\em twisted square} in $h$ if there exist indexes $i,j \in \{0, \ldots, 2^n-1\}$ such that 
$$h_i h_{i+1} h_j h_{j+1} = v_0v_1v_3v_2.$$
\item We will say that the edges $\edg{v_0}{v_1}$ and $\edg{v_2}{v_3}$ are the {\em rims} and the edges $\edg{v_1}{v_2}$ and $\edg{v_3}{v_0}$ are the {\em rays}.
\end{enumerate}
In figure~\ref{fg.02} we illustrate those notions, the rays are displayed as dashed lines, the rims as continuous lines since they actually ``are part of'' the Hamiltonian cycle.
We will say that the Hamiltonian cycle $h$ {\em contains a square}, or that a square is {\em inscribed} in $h$, if a straight or twisted square is contained in $h$. A  {\em square free Hamiltonian} cycle contains no squares.

\begin{conjecture}[Inscribed squares in Hamiltonian cycles.] \label{c:squaresConjecture} Let $n \ge 2$ be an integer. No Hamiltonian cycle of $\hyp{n}$ is square free, or equivalently any Hamiltonian cycle has inscribed a square.
\end{conjecture}

Since any 4-cycle in $\hyp{n}$ contains two alternating dimensions, the rims (or rays) in each  square inscribed in a Hamiltonian cycle of $\hyp{n}$ belong to the same dimension. 

Now let us establish some conditions in order to have edges of the same dimension in a square contained within a Hamiltonian cycle of $\hyp{n}$.

From the construction of graph ${\cal D}_{in}$, specially the edge definition at relation~(\ref{eq.bied}), and the Theorem~\ref{p:isomorphism}, if the chromatic vector of some Hamiltonian cycle $h$ of $\hyp{n}$ contains entries greater than the independence number of $\hyp{n-1}$ then $h$ must contain a square. Since $\alpha(\hyp{n-1}) = 2^{n-2}$ the following theorem results.

\begin{theorem} \label{p:independence}
Let $n\ge 2$ be an integer, and let $h$ be a Hamiltonian cycle of $\hyp{n}$. If some entry $i$ in the chromatic vector $c(h)$ is greater than $2^{n-2}$ then $h$ contains a square whose rims are \idim\  edges. 
\end{theorem}

Theorem~\ref{p:independence} is a general result aiming to prove Conjecture~\ref{c:squaresConjecture} for Hamiltonian cycles whose chromatic vectors have entries bigger than $2^{n-2}$. Let us reduce further this lower bound.

Let $G$ be a bipartite graph. The {\em equi-independence} number of $G$, denoted $\alpha_=(G)$, is the cardinality of the maximum independent set in $G$ having half of its vertexes in one bipartition class of $G$ and half in the other. 

According to the proofs of Theorem~\ref{p:independence} and Theorem~\ref{p:bipartition} we have:

\begin{corollary} \label{c:equiIndependence}
Let $n \ge 2$ be an integer, and let $h$ be a Hamiltonian cycle of $\hyp{n}$. If some entry $c_i$ in $c(h)$ is greater than $\alpha_=(Q_{n-1})$ then $h$ contains a square whose rims are \idim\  edges. 
\end{corollary}

Let us estimate a lower bound for the equi-independence number of the $n$-dimensional hypercube.

\begin{theorem} \label{p:lowerBound}
Let $n \ge 3$ be an integer. Then $\alpha_=(\hyp{n}) \ge 2^{n-2}$.
\end{theorem}

\noindent{\em Proof}: For any two bits $b_0,b_1 \in \{0,1\}$ let $\hyp{nb_0b_1}$ be the subgraph of $\hyp{n}$ induced by the vertexes in $V(\hyp{n})$ having values $b_0$ and $b_1$ at their first two coordinates. The resulting four graphs are isomorphic to $\hyp{n-2}$. 

Let $b\in\{0,1\}$, and let $I_b = V(Q_{nbb}) \cap (b\mbox{-bipartition class of } \hyp{n})$. No vertex in $I_0$ is adjacent to a vertex of $I_1$ and in fact $I = I_0 \cup I_1$ is an independent set of $\hyp{n}$. Moreover, each vertex in $V(\hyp{nbb})$ not contained in $I$ is adjacent to one vertex in $\hyp{nbb}$, and each vertex in $\hyp{n01}$ (or in $\hyp{n10}$) is adjacent to a vertex in $I_0$ (or a vertex in $I_1$.) 

In order to illustrate the last claim let $v = (0,1,b_2, \ldots, b_{n-1})$ be an arbitrary vertex in $\hyp{n,0,1}$, if $\mbox{\it par}_n(b_2, \ldots, b_{n-1}) = 0$ then $v$ is adjacent to $(0,0,b_2, \ldots, b_{n-1})$ which is in $I_0$, otherwise $v$ is adjacent to $(1,1,b_2, \ldots, b_{n-1})$ which is in $I_1$.

The final conclusion is that $I$ is a maximal independent set in $\hyp{n}$. Since $\mbox{card}(I_0) = \mbox{card}(I_1) = 2^{n-3}$ and $$\mbox{card}(I) = \mbox{card}(I_0)+\mbox{card}(I_1) = 2^{n-3}+2^{n-3} = 2^{n-2}$$ the Theorem follows. \hfill $\smallbox$

Now let us check that finding the equi-independence number of a bipartite graph $G$ is not more difficult than evaluating the independence number of a graph $G'$ derived from $G$.

\begin{theorem} \label{p:equiIndBipartite}
The problem of finding the equi-independence number of a bipartite graph is polynomially-time reducible to the problem of finding the independence number of a graph. 
\end{theorem}

\noindent{\em Proof}: Let $G$ be a bipartite graph, and let $V_0$, $V_1$ be the bipartition classes of $G$. The graph $G' = (V',E')$, with set of vertexes 
$V'=\{(v_0,v_1) \in V_0 \times V_1 | (v_0, v_1) \notin E(G)\}$ and set of edges $$E' = \{((v_0, v_1), (v_0', v_1')) \in V(G')\times V(G') | (v_0, v_1') \in E(G) \mbox{ or } (v_1, v_0') \in E(G) \mbox{ or } v_0 = v_0' \mbox { or } v_1 = v_1'\}),$$ has an independent set $I'$ if and only if the set 
$I = \{v \in V(G) | v \mbox{ is an element in some vertex of } I'\}$ is an equi-independent set of $G$. In this construction $\mbox{card}(I) = 2\,\mbox{card}(I')$. Since $G'$ can be built in polynomial time, the proposition follows. \hfill$\smallbox$. 

Theorem~\ref{p:equiIndBipartite} was used to compute the equi-independence number of the hypercubes of dimension three to seven. The results are summarized in Table~\ref{t:results}. 
Vertexes in the maximal equi-independence sets in Table~\ref{t:results} are coded in binary. The equi-independence numbers for $n$ equals to four and six reach the lower bound in Theorem~\ref{p:lowerBound} and so this bound is tight. The graph $\hyp{n}'$ grows very fast as $n$ is increased, thus for $n = 8$ the corresponding value is not included. 

\begin{table}[ht]
\caption{Equi-independence number of $\hyp{n}$ for small values of $n$.}
\label{t:results}\centering 
\begin{tabular}{ |c|r|l|r|r|}
\hline
  n & $\alpha_=(\hyp{n})$ & A maximal equi-independent set & $|V(\hyp{n}')|$ & $|E(\hyp{n}')|$ \\
  \hline
  3 & 2 & \{0, 7\}        & 4 & 6 \\
  4 & 4 & \{0, 7, 9, 14\} & 32 & 448 \\
  5 & 10 & \{0, 7, 9, 19, 10, 21, 12, 22, 24, 31\}                         & 176  & 9720 \\
  6 & 16 & \{0, 7, 9, 19, 33, 21, 10, 22, 34, 28,                          & 882  & 137536\\
    &    & \hspace{0.03in} 36, 56, 43, 31, 45, 55\}                        &      &\\
  7 & 40 & \{0, 7, 9, 19, 33, 67, 10, 21, 34, 69,                          & 3648 & 1577184\\
    &    & \hspace{0.03in} 12, 81, 36, 22, 24, 70, 40, 82, 72, 84,         &      &\\ 
    &    & \hspace{0.03in} 48, 31, 96, 47, 57, 79, 105, 55, 58, 87,        &      &\\
    &    & \hspace{0.03in} 106, 103, 60, 91, 108, 115, 120, 93, 117, 126\} &      &\\
\hline  
\end{tabular}
\end{table}

The study on the inscribed square conjecture in the Hamiltonian cycles of the hypercube concludes as follows.

\begin{theorem}
Conjecture~\ref{c:squaresConjecture} holds for $2 \le n \le 7$.
\end{theorem}

\noindent {\em Proof}: From Table~\ref{t:results} we know that $\alpha_=(\hyp{6}) = 16$. It is impossible that the chromatic vector of a Hamiltonian cycle of $\hyp{7}$ have all its entries lower or equal that 16; otherwise the sum of the seven entries would be at most $16\cdot 7 = 112$, but we know from Lemma~\ref{l:basicChromatic} that this sum should be 128. Thus some entry should be greater than 16. It follows from Corollary~\ref{c:equiIndependence} that Conjecture~\ref{c:squaresConjecture} is true for $n=7$. The proof is analogous for values of $n$ lower than $7$. $\smallbox$

Computer experimentation suggests that the inscribed square conjecture in Hamiltonian cycles is true for any value of $n$.

\section{Conclusions}

Theorem~\ref{p:bipartition} entails several applications, e.g., in~\cite{Fink07} it is posed the following question: The partial matchings in $\hyp{n}$ can be extended to Hamiltonian cycles? We are able to submit a partial answer to this question: if the partial matching contains edges in the same dimension that violate the equilibrium condition of the Theorem~\ref{p:bipartition}, and it is not allowed to incorporate more edges in the same dimension, then no such extension exists. Other application of the Theorem~\ref{p:bipartition} consists in a pruning strategy to generate exhaustively all the Hamiltonian cycles of $\hyp{n}$.

On the other hand several problems and conjectures have remained open, and their study could reveal important structural properties of $\hyp{n}$. In the following presentation we will refer to the operators introduced at the end of section~\ref{s:coloring}.

As conjectures, besides the already stated Conjecture~\ref{c:squaresConjecture}, let us state the following:
\begin{itemize}
\item {\em Let $n\geq 2$ and let $k$ be an even integer with $2 \le k \le 2^{n-1}$. A vector $v \in \hyp{k}$ is the parity vector of a Hamiltonian cycle of $\hyp{n}$, i.e. $v=\beta_i(h)$ for some $i$ and Hamiltonian cycle $h$, if and only if it contains the same number of $0$'s and $1$'s.}

In the Theorem~\ref{p:bipartition}, the ``only if'' part of this conjecture was proved.
\item {\em The lower bound in Theorem~\ref{p:lowerBound} is reached for all even values of $n$.}
\end{itemize}

As open problems, let us state the following:
\begin{itemize}
\item {\em Characterize the chromatic vectors of Hamiltonian cycles of the hypercube.} Or, in other words, {\em characterize the image of the map $c:\{\mbox{Hamiltonian cycles}\}\to\N^n$.}

Lemma~\ref{l:basicChromatic} states some conditions that are just necessary. 

\item {\em Characterize the chromatic-segments vectors of Hamiltonian cycles of $\hyp{n}$. } Or, in other words, {\em characterize, for each $i$, the image of the map $\delta_i:\{\mbox{Hamiltonian cycles}\}\to\N^n$.}

A positive answer to the first of the above stated conjectures will give immediately necessary conditions for a vector to be a chromatic-segments vector of a Hamiltonian cycle of the hypercube. But it would be far from a full characterization.
\item {\em Calculate $\alpha_=(\hyp{n})$.}
\end{itemize}

\end{document}